\documentclass[12pt,reqno]{amsart}
\usepackage{amsaddr}
\usepackage{hyperref}
\hypersetup{colorlinks=true, linkcolor=blue, citecolor=red, urlcolor=blue}
\usepackage[capitalize,nameinlink]{cleveref}
\author{Tirthankar Bhattacharyya, Ritul Duhan and Chandan Pradhan}
\address{Department of Mathematics, 	Indian Institute of Science, 	Bangalore 560012, India}
\email{tirtha@iisc.ac.in, ritul2023@iisc.ac.in, chandan.pradhan2108@gmail.com}

\usepackage{color, bm, amscd, tikz-cd}
\usepackage{csquotes}
\usepackage{tikz}

\newtheorem{thm}{Theorem}[section]

\newtheorem{cor}[thm]{Corollary}
\newtheorem{lem}[thm]{Lemma}

\newtheorem{obs}[thm]{Observation}

\newtheorem{prop}[thm]{Proposition}

\newtheorem{theorem}{Theorem}
\newtheorem{defn}[thm]{Definition}
\newtheorem{rem}[thm]{Remark}

\numberwithin{equation}{section}

\def\tilde{\widetilde}

\def\B{\mathcal{B}}
\def\bar{\overline}
\newcommand{\bra}[1]{\langle #1 |}
\newcommand{\ket}[1]{| #1 \rangle}
\def\C{\mathbb{C}}

\def\E{\mathcal{E}}

\def\epsilon{\varepsilon}

\def\H{\mathcal{H}}

\def\K{\mathcal{K}}

\def\N{\mathbb{N}}

\def\phi{\varphi}

\def\cT{\mathcal{T}}

\newcommand{\UE}[2]{\mathcal{UCP}_{#2}\big(#1\big)}

\def\Z{\mathbb{Z}}

\newcommand{\lr}[2]{\left\langle #1, #2\right\rangle}

\def\spb{\textit{spb}}

\begin{document}
	
	\title[Gromov-Hausdorff convergence of metric spaces of UCP maps]{Gromov-Hausdorff convergence of metric spaces of UCP maps}
	
	\maketitle
	\begin{abstract}
		It is shown that van Suijlekom's technique of imposing a set of conditions on operator system spectral triples ensures Gromov-Hausdorff convergence of sequences of sets of unital completely positive maps (equipped with the BW-topology which is metrizable). This implies that even when only a part of the spectrum of the Dirac operator is available together with a certain truncation of the $C^*$-algebra, information about the geometry can be extracted.  
	\end{abstract}
	
	{\footnotesize \noindent 2020 Mathematics Subject Classification:  46L07, 46L87, 47C15, 47L25, 58B34. \\
		Keywords: Operator system, UCP maps,  Spectral truncation, Gromov--Hausdorff convergence

	\section{Introduction}\label{sec:intro}
Let $C(S^1)$ be the commutative $C^*$-algebra of all continuous functions on the unit circle $S^1$. The Fejer-Riesz operator system $C(S^1)_{(n)}$ for any $n \ge 1$ is the vector space of all continuous functions $f$ on the unit circle $S^1$ such that the Fourier coefficients $\hat f (k)$ vanish for all $k \in \mathbb Z$ with $|k| \ge n$, i.e., $f \in C(S^1)_{(n)}$ if and only if 
$$ f(z) = \sum_{k=-(n-1)}^{n-1} a_k z^k$$
for $z \in S^1$. As a subspace of the $C^*$-algebra $C(S^1)$, it is $*$-closed because the involution is given by
$$ \sum_{k=-(n-1)}^{n-1} a_k z^k \rightarrow \sum_{k=-(n-1)}^{n-1} \bar{a_{-k}} z^k$$
and contains the identity of the $C^*$-algebra. So, it is an operator system.

On the other hand, the $n \times n$ Toeplitz matrices over $\mathbb C$, denoted by $\mathcal T_n$ form an operator system because it is a $*$-closed unital subspace of the $C^*$-algebra of all $n \times n$ matrices over $\mathbb C$. 

Consider the linear map $\delta : \mathcal T_n \rightarrow C(S^1)_{(n)}$ which takes a Toeplitz matrix $t = (( \tau_{k-l} ))_{k,l=0}^{n-1}$ to the linear functional $\varphi_t$ on $C(S^1)_{(n)}$ defined by
$$ \varphi_t(f) = \sum_{k=-(n-1)}^{n-1} \tau_{-k} a_k$$
where $ f(z) = \sum_{k=-(n-1)}^{n-1} a_k z^k$. Connes and van Suijlekom in Proposition 4.6 of \cite{CS20} showed that this is a linear unital order isomorphism between the dual of $C(S^1)_{(n)}$ and $\mathcal T_n $. In Theorem 1.1 in \cite{Far}, Farenick showed that this duality identification is unital completely isometrically isomorphic.

We are motivated by the convergence results in the literature, e.g., the Gromov-Hausdorff convergence of quantum metric spaces in the seminal paper \cite{Rie00}. More recently, in \cite{CS20} and  \cite{Suij-JGP}, state spaces of the operator systems $C(S^1)_{(n)}$ and $\mathcal T_n$ were dealt with. It was shown that in spite of the duality mentioned above, the sequences of state spaces $\mathcal S (C(S^1)_{(n)})$ and $\mathcal S (\mathcal T_n)$ both converge to the state space $\mathcal S (C(S^1))$. More generally, in \cite{Suij-JGP}, general conditions were found on operator system spectral triples that guarantee the  Gromov-Hausdorff convergence of state spaces.

Given a closed operator system $E$ and a (complex, separable) Hilbert space $\K$, let $\mathcal{UCP}_{\K}(E)$ denote the collection of all unital completely positive (UCP) maps from $E$ into the algebra of all bounded operators $\B(\K)$. We shall show in this note that the general conditions due to van Suijlekom mentioned above actually guarantee much more, viz., the  Gromov-Hausdorff convergence of sequences of sets of UCP maps. In particular, $\mathcal{UCP}_{\K}(C(S^1)_{(n)})$ and $\mathcal{UCP}_{\K}(\mathcal T_n)$ both converge in the Gromov-Hausdorff metric to $\mathcal{UCP}_{\K}(C(S^1))$ as $n$ tends to infinity, implying that for large $n$, $\mathcal{UCP}_{\K}(C(S^1)_{(n)})$ and $\mathcal{UCP}_{\K}(\mathcal T_n)$ are close in the Gromov-Hausdorff metric. 

The main result is Theorem \ref{thm:UCP-GH-conv} where we consider a sequence of  operator system spectral triples $\{ (\E_n, \H_n, D_n) \}_n$ related to  an operator system spectral triple $(\E, \H, D)$ by a $\mathrm{C}^1$-approximate complete order isomorphism. Section \ref{appl} contains applications of the theorem to examples like the ones mentioned above as well as for {\em polyhedral truncations} motivated by works of Travaglini \cite{Trav}. The $\{ (\E_n, \H_n, D_n) \}_n$ are often called spectral truncations of $(\E, \H, D)$.  Spectral truncations have a widespread appeal - see  \cite{FLP}, \cite{LeSu24}, \cite{Suij-JGP} and \cite{Rie2023} for noncommutative geometry, \cite{Leim} for compact quantum groups, and \cite{HHIK} for machine learning.

	\section{Sets of CP maps as metric spaces}
\subsection{The BW-topology}
	To talk of the Gromov-Hausdorff metric, we need to make a metric space out of $\mathcal{UCP}_{\K}(E)$ first. For the following, a good reference is \cite{Pau02}. 
	
	\begin{defn}[BW-topology]
		For Banach spaces $X$ and $Y$, let $\mathcal{B}(X,Y^*)$ denote the bounded linear maps of $X$ into $Y^*$. For fixed vectors $x$ in $X$ and $y$ in $Y$, define a linear functional $x \otimes y$ on $\mathcal{B}(X,Y^*)$ by $x \otimes y(L)$ = $L(x)(y)$. Let $Z$ denote the closed linear span in $(\mathcal{B}(X,Y^*))^*$ of these elementary tensors. Then $\mathcal{B}(X,Y^*)$ is isometrically isomorphic to $Z^*$ with the duality given by \\
		\[ \langle L, x \otimes y \rangle = x \otimes y (L)   \] 
		The weak$^*$ topology that is induced on $\mathcal{B}(X,Y^*)$ by this identification is called the BW topology $($for bounded weak$)$. 
	\end{defn}

If  $Y^*$ is the algebra $\mathcal{B}(\mathcal{H})$ of bounded operators on a Hilbert space $\mathcal{H}$, then a bounded net $L_{\lambda}$ in $\mathcal{B}(X,\mathcal{B}(H))$ converges in the BW topology to $L$ if and only if $\langle L_{\lambda}(x)h,h^\prime \rangle$ converges to $\langle L(x)h,h^\prime \rangle$ for all $h$, $h^\prime$ in $\mathcal{H}$ and $x$ in $X$. Naturally, $\mathcal{UCP}_{\K}(E)$ inherits the BW-topology from $\mathcal{B}(E,\mathcal{B}(K))$. It is known that $\mathcal{UCP}_{\K}(E)$ is compact in the BW-topology \cite[Theorem 7.4]{Pau02}.

\subsection{A metric for the BW-topology when $E = C(\Omega)$.} 

Let $(\Omega, d_{\Omega})$ be a compact metric space. Consider the set of all Lipschitz functions 
\begin{align} \label{def:Lip}
\text{Lip}(\Omega) := \{f: \Omega \to \C \text{ such that } \|f\|_1 := \sup_{x,y \in \Omega} \frac{|f(x) - f(y)|}{d_{\Omega}(x, y)} < \infty\}.
\end{align}
Let $\{k_n \in \K : n \in \N\}$ be a fixed dense subset of the closed unit ball of $\K$. Define $d: \mathcal{UCP}_{\K}(C(\Omega)) \times \mathcal{UCP}_{\K}(C(\Omega)) \to [0, \infty)$ by

	\begin{align}\label{def:metric}
		d(\phi, \psi) = \sup\left\{\sum_{p=1}^{\infty} \frac{1}{2^p} \left|\langle (\phi - \psi)(f) k_p, k_p \rangle \right| : f \in \text{Lip}(\Omega), \|f\| \leq 1, \|f\|_1 \leq 1 \right\}.
	\end{align}
	\begin{lem}
		The function $d$ given by \eqref{def:metric} forms a metric on $\mathcal{UCP}_{\K}(C(\Omega))$.
	\end{lem}
	
	\proof
	Clearly, symmetry and the triangle inequality hold. Also, the implication $\phi = \psi \implies d(\phi, \psi) = 0$ is immediate. Suppose now that $d(\phi, \psi) = 0$. Then for all $f \in \text{Lip}(\Omega)$ and for all $p$, we have $\lr{(\phi - \psi)(f) k_p}{k_p} = 0$. By denseness, this implies that 
	\begin{align}\label{positivity1}
		 \lr{(\phi - \psi)(f) k}{k} = 0 \text{ for all } k \in \K. 
	\end{align}
	For $h, k \in \K$, \eqref{positivity1} and the polarization identity
	\begin{align*}
		\lr{(\phi - \psi)(f) h}{k} = \sum_{l=0}^{3} \frac{i^l}{4} \lr{(\phi - \psi)(f)(h + i^l k)}{h + i^l k}
	\end{align*}
	together yield $\lr{(\phi - \psi)(f) h}{k} = 0$. Hence $(\phi - \psi)(f) = 0$ for all $f \in \text{Lip}(\Omega)$, and consequently $\phi = \psi$. 
	This completes the proof. 
	\endproof

The following is probably known. We give a proof because we could not find a suitable reference in the literature.

	\begin{lem}\label{lem:BW=d}
	The BW-topology on $\mathcal{UCP}_{\K}(C(\Omega))$ coincides with the metric topology on $\mathcal{UCP}_{\K}(C(\Omega))$ induced by the metric $d$ defined in \eqref{def:metric}.
	\end{lem}
	
	\proof
	Let $(\mathcal{UCP}_{\K}(C(\Omega)), d)$ and $(\mathcal{UCP}_{\K}(C(\Omega)), BW)$ denote the metric topology and the {\it BW} topology on $\UE{E}{K}$, respectively. To show that these topologies are the same, it is enough to show that the identity map $I: (\mathcal{UCP}_{\K}(C(\Omega)), BW) \to (\mathcal{UCP}_{\K}(C(\Omega)), d)$ is a homeomorphism. However, since $(\mathcal{UCP}_{\K}(C(\Omega)), BW)$ is a compact space and $(\mathcal{UCP}_{\K}(C(\Omega)), d)$ is Hausdorff, it suffices to show the continuity of $I$.
	
	Let $\{\phi_\lambda\}_{\lambda \in \Lambda}$ be a net in $(\mathcal{UCP}_{\K}(C(\Omega)), BW)$ such that $\phi_\lambda$ converges to some $\phi$ in $(\mathcal{UCP}_{\K}(C(\Omega)), BW)$. Our claim is that $\phi_\lambda$ also converges to $\phi$ in $(\mathcal{UCP}_{\K}(C(\Omega)), d)$. Suppose not. Then there exist $\epsilon > 0$ and an infinite subset $\hat{\Lambda} \subseteq \Lambda$ such that $d(\phi_\lambda, \phi) > 2\epsilon$ for all $\lambda \in \hat{\Lambda}$. Invoking the Axiom of Choice, we may assume that $\hat{\Lambda} = \mathbb{N}$. Thus, we have $d(\phi_\lambda, \phi) > 2\epsilon$ for $\lambda \in \mathbb{N}$. Therefore, there exists a sequence $\{f_\lambda\}$ in $\text{Lip}(X)$ with $\|f_\lambda\| \leq 1$ and $\|f_\lambda\|_1 \leq 1$ such that
	
	\begin{align}\label{bw:1}
		\sum_{p=1}^{\infty} \frac{1}{2^p} \left|\langle (\phi_\lambda - \phi)(f_\lambda) k_p, k_p \rangle \right|>2\epsilon.
	\end{align} 
	Since $\|f_\lambda\|_{ 1}\leq 1$, so $\{f_\lambda\}$ is a sequence of equicontinuous functions. Next, in view of 
	\[\left|\langle (\phi_\lambda - \phi)(f_\lambda) k_p, k_p \rangle \right|\leq 2\]
	 and the fact that for a given $\epsilon>0$, there exist an $N_0\in\N$ such that 
	\[\sum_{p=N_0+1}^{\infty}\frac{1}{2^{p-1}}<\epsilon,\]
	we have 
	\begin{align}\label{bw:2}
		\sum_{p=N_0+1}^{\infty} \frac{1}{2^p} \left|\lr{ (\phi_\lambda - \phi)(f_\lambda) k_p }{ k_p}\right|  <\epsilon
	\end{align} 
for every $\lambda$. Therefore, \eqref{bw:1} and \eqref{bw:2} together give
	\begin{align}\label{bw:3}
		\sum_{p=1}^{N_0} \frac{1}{2^p} \left|\lr{ (\phi_\lambda - \phi)(f_\lambda) k_p }{ k_p}\right| >\epsilon \;\; \text{ for every $\lambda$}.
	\end{align}
	Since there are only finitely many $p$ in \eqref{bw:3}, we can choose a $p_0 \in \{1, \ldots, N_0\}$ such that after passing to a subsequence and continuing to denote the subsequence by $\lambda$, we have
	\begin{align}\label{bw:4}
		\left|\lr{ (\phi_\lambda - \phi)(f_\lambda) k_{p_0} }{ k_{p_0}}\right| >\frac{\epsilon}{N_0}.
	\end{align}
	Since $\{f_\lambda\}$ is a uniformly bounded sequence in $E$ of equicontinuous functions, using the Arzel\`a-Ascoli theorem there exists a convergent subsequence (which we still denote by $\{f_\lambda\}$) that converges to some $f \in E$. Therefore, there exists an $M \in \mathbb{N}$ such that
	\begin{align}\label{bw:5}
		\|f_\lambda-f\|<\frac{\epsilon}{4N_0} \text{ for all } \lambda\geq M.
	\end{align}
	Finally, by \eqref{bw:3}, \eqref{bw:4}, and \eqref{bw:5}, we have
	\begin{align*}
		&\left|\lr{ (\phi_\lambda - \phi)(f) k_{p_0} }{ k_{p_0}}\right|\\
		=& \left|\left(\lr{ (\phi_\lambda - \phi)(f_\lambda) k_{p_0} }{ k_{p_0}}\right)+ \left(\lr{\phi(f_\lambda-f) k_{p_0} }{ k_{p_0}}+ \lr{ \phi_\lambda(f-f_\lambda) k_{p_0} }{ k_{p_0}}\right)\right|\\
		\geq& \Big|\left|\left(\lr{ (\phi_\lambda - \phi)(f_\lambda) k_{p_0} }{ k_{p_0}}\right)\right|- \left|\left(\lr{\phi(f_\lambda-f) k_{p_0} }{ k_{p_0}}+ \lr{ \phi_\lambda(f-f_\lambda) k_{p_0} }{ k_{p_0}}\right)\right|\Big|\\
		>& \frac{\epsilon}{N_0}-\frac{2\epsilon}{4N_0} = \frac{\epsilon}{2N_0},
	\end{align*}
	which contradicts the fact that $\{\phi_\lambda\}$ converges to $\phi$ in $(\mathcal{UCP}_{\K}(C(\Omega)), BW)$. Hence, $I$ is continuous. This completes the proof.
		\endproof

\subsection{A metric for the BW-topology when $E$ is a finite-dimensional operator system}

The metric in \eqref{def:metric} depended on the class in \eqref{def:Lip}. In an operator system, we do not have access to such a class in general. To compensate for that loss, we need a {\em Dirac operator} $D$. 

	\begin{defn}An operator system spectral triple is a triple $(\E, \H, D)$ where $\E$ is a dense subspace of an operator system $E$ in $\mathcal{B}(\H)$, $\H$ is a Hilbert space, and $D$ is a self-adjoint operator in $\H$ with compact resolvent, such that $[D, T]$ is a bounded operator for all $T \in \E$.
	\end{defn}
It is well-known that the Dirac operator induces a {\em Lipschitz semi-norm}: $\| T\|_1 = \| [D,T]\|$ on $\E$. Hence the analogue of \eqref{def:metric} in this case is 
\begin{align*}
		d_E(\phi, \psi) = \sup\left\{\sum_{p=1}^{\infty} \frac{1}{2^p} \left|\langle (\phi - \psi)(f) k_p, k_p \rangle \right| : f \in \E, \|f\| \leq 1, \|f\|_1 \leq 1 \right\}
	\end{align*}
where $\phi, \psi$ are in $\mathcal{UCP}_{\K}(E)$. 

	\begin{defn}
		We say that an operator system spectral triple $(\E, \H, D)$ has the {\it BW} property if the metric topology on $\UE{E}{\K}$ induced by the metric $d_{E}$ coincides with the {\it BW}-topology on $\UE{E}{\K}$.
	\end{defn}
	\begin{lem}\label{lem:BW=d:1}
		Let $(\E, \H, D)$ be an operator system spectral triple such that $\E$ is finite-dimensional. Then $(\E, \H, D)$ has the {\it BW} property.
	\end{lem}
	\begin{proof}
		We are omitting the proof because it is similar to the proof of Lemma \ref{lem:BW=d}. The compactness argument which involved the  Arzel\`a-Ascoli theorem in the proof of Lemma \ref{lem:BW=d} is automatic in the finite-dimensional case because of compactness of the closed unit ball.
	\end{proof}

\begin{obs}
Obviously, in the finite dimensional case, one can omit the condition $\|f\|_1 \leq 1$ in the definition of the metric. The metric could be defined as 
$$d_E^\prime(\phi, \psi) = \sup\left\{\sum_{p=1}^{\infty} \frac{1}{2^p} \left|\langle (\phi - \psi)(f) k_p, k_p \rangle \right| : f \in \E, \|f\| \leq 1 \right\}.$$ 
This too would metrize the BW-topology. So, the two metrics $d_E$ and $d_E^\prime$ are equivalent. In other words, the Dirac operator is not necessary to define the metric. We introduce the Dirac operator at this stage because it is essential in formulating the sufficient conditions which will guarantee convergence in the next section. 
\end{obs}
	
	\section{The Gromov--Hausdorff convergence}

In this section, we shall consider the usual {\em Gromov-Hausdorff} metric $d_{GH}(M_1, M_2)$ between two compact metric spaces $(M_1, d_1)$ and $(M_2, d_2)$ and apply it to measure distances between operator system spectral triples. For various notions of distances between Lip-normed operator systems and inter-relations between them, we refer the reader to section 3 of \cite{KL}. 

\begin{defn}$d_{GH}(M_1, M_2)$ is the infimum of the set
\begin{align*}
\{ & d_H(f(M_1), g(M_2)) : M \text{ is a compact metric space and } \\
&f \text{ and } g \text{ are isometric embeddings of } M_1 \text{ and } M_2 \text{ into } M\}\end{align*}
where $d_H$ stands for the Hausdorff distance between closed subsets of a compact metric space. 
\end{defn}

For proving convergence in the Gromov-Hausdorff metric, we use the {\em distortion}.
	
	\begin{defn}\cite[Definition~7.3.1]{BBI01}\label{correspondence}
		Let $M_1$ and $M_2$ be two sets. A total onto correspondence between $M_1$ and $M_2$ is a relation $\mathfrak{R} \subseteq M_1 \times M_2$ such that for every $x \in M_1$ there exists at least one $y \in M_2$ with $(x,y) \in \mathfrak{R}$ and similarly for every $y \in M_2$ there exists an $x \in M_1$ with $(x,y) \in \mathfrak{R}$.
	\end{defn}

	\begin{defn}\cite[Definition~7.3.10]{BBI01}
		Let $\mathfrak{R}$ be a total onto correspondence between metric spaces $(M_1, d_1)$ and $(M_2, d_2)$. The \textit{distortion} of $\mathfrak{R}$ is defined by
		\[\operatorname{dis } \mathfrak{R} = \sup \{|{d_1(x, x') - d_2(y, y')|} : (x,y), (x', y') \in \mathfrak{R} \}.\]
	\end{defn}

	We depend on Theorem 7.3.25 of \cite{BBI01} which gives the crucial relation that for any two metric spaces $(M_1, d_1)$ and $(M_2, d_2)$, 
\begin{align} \label{eqn:dist-GH} d_{GH}(M_1,M_2) = \frac{1}{2} \inf_{\mathfrak{R}}(\operatorname{dis } \mathfrak{R}),\end{align}
where the infimum is taken over all total onto correspondences $\mathfrak{R}$ between $M_1$ and $M_2$.

	\begin{defn}
		Let $\{ (\E_n, \H_n, D_n) \}_n$ be a sequence of operator system spectral triples and let $(\E, \H, D)$ be an operator system spectral triple. A $\mathrm{C}^1$-approximate complete order isomorphism is a pair  $(R_n, S_n)$ of linear maps $R_n: \E \to \E_n$ and $S_n: \E_n \to \E$ for each $n$, satisfying the following conditions:
		\begin{enumerate}
			\item The maps $R_n$ and $S_n$ are unital completely positive (ucp) maps.
			\item $R_n$ and $S_n$ are \emph{$\mathrm{C}^1$-contractive}, i.e.\@ they are contractive both with respect to norm $\|\cdot\|$ and with respect to Lipschitz seminorm $\|[D,\cdot]\|$, respectively $\|[D_n,\cdot]\|$.
			\item There exist sequences $c_n$ and $c_n'$ both converging to zero such that
			\begin{align*}
				\| S_n \circ R_n (f) - f \| &\leq c_n \| f \|_1, \\
				\| R_n \circ S_n (h) - h \| &\leq c_n' \| h \|_1
			\end{align*}
			for all $f \in \E$ and $h \in \E_n$.
		\end{enumerate}
	\end{defn}
	
	Since the composition of two ucp maps is again a $\mathrm{ucp}$ map, we can pull back ucp maps as follows:
	
	\begin{align*}
		&R_n^* :  \UE{E_n}{\K} \to \mathcal{UCP}_{\K}(E); \qquad \phi_n \mapsto  \phi_n \circ R_n,\\
		&S_n^* :  \UE{E}{\K} \to \mathcal{UCP}_{\K}(E_n); \qquad \phi\mapsto \phi \circ S_n.
	\end{align*}

	\smallskip
	
	The following results transfer \cite[Proposition 4]{Suij-JGP} to the setting of $\mathrm{ucp}$ maps.
	
	\begin{prop}\label{prop:estimate}
		If $(R_n,S_n)$ is a $C^1$-approximate complete order isomorphism for $(\E_n,\H_n,D_n)$ and $(\E,\H,D)$, then 
		\begin{enumerate}
			\item\label{inequality: first} For all $\phi,\psi \in \mathcal{UCP}_{\K}(E)$ we have
			\[
			d_{E_n}(\phi \circ S_n, \psi \circ S_n) \leq d_E (\phi,\psi) \leq d_{E_n}(\phi \circ S_n, \psi \circ S_n) +2 c_n.
			\]
			
			\item\label{inequality: second} For all $\phi_n,\psi_n \in \mathcal{UCP}_{\K}(E_n)$ we have
			\[d_E(\phi_n \circ R_n, \psi_n \circ R_n) \leq d_{E_n} (\phi_n,\psi_n) \leq d_E(\phi_n\circ R_n, \psi_n \circ R_n) +2 c_n'.\]
			
		\end{enumerate}
		
	\end{prop}

	\proof
	Let $\phi, \psi \in \mathcal{UCP}_{\K}(E)$. Note that, for $n \in \N$, $S_n: E_n \to E$ are Lipschitz contractive as well as contractive; therefore, we have:

	\begin{align}\label{inequality:1}
		\nonumber&d_{E_n}(\phi \circ S_n, \psi \circ S_n)\\
		\nonumber=&\sup\left\{\sum_{p=1}^{\infty} \frac{1}{2^p} \left|\langle (\phi - \psi)(S_nf) k_p, k_p \rangle \right| : f \in \E_n, \|f\| \leq 1, \|f\|_1 \leq 1 \right\}\\
		\leq & \sup\left\{\sum_{p=1}^{\infty} \frac{1}{2^p} \left|\langle (\phi - \psi)(g) k_p, k_p \rangle \right| : g \in \E, \|g\| \leq 1, \|g\|_1 \leq 1 \right\} = d_{E}(\phi, \psi)
	\end{align}
	On the other hand, for all $f \in \E$ with $\|f\|,\, \|f\|_1 \leq 1$, we have
	\begin{align}\label{inequality:2}
		\nonumber&\sum_{p=1}^{\infty} \frac{1}{2^p} \left|\langle (\phi - \psi)(f) k_p, k_p \rangle \right|\\
		\nonumber&\leq \sum_{p=1}^{\infty} \frac{1}{2^p}|\phi(S_n(R_n(f))) - \psi(S_n(R_n(f))) |  + \sum_{p=1}^{\infty} \frac{1}{2^p}| \phi(f)- \phi(S_n(R_n(f)))|\\
		\nonumber& + \sum_{p=1}^{\infty} \frac{1}{2^p}| \psi(f) - \psi(S_n(R_n(f)))| \\
& \leq  d_E(\phi\circ S_n, \psi \circ S_n) + 2 c_n,
	\end{align}
	as $\| \phi\| = \|\psi\|=1$ and $\|R_n(f)\|_1 \leq \| f \|_1 \leq 1$. Therefore, \eqref{inequality:1} and \eqref{inequality:2} together yield \eqref{inequality: first}. By similar lines of argument, we also prove \eqref{inequality: second}.
	\endproof

	Now we state and prove our main result below.
	\begin{theorem}\label{thm:UCP-GH-conv}
		Let $n\in\N$. Let $(\E_n, \H_n, D_n)$, and $(\E, \H, D)$ be operator system spectral triples having \textit{BW}-property. If $(R_n, S_n)$ is a $C^1$-approximate complete order isomorphism for $(\E_n, \H_n, D_n)$ and $(\E, \H, D)$, then $\mathcal{UCP}_{\K}(E_n)$ converges to $\mathcal{UCP}_{\K}(E)$ in the Gromov--Hausdorff metric.
	\end{theorem}
	
	\proof
	
	For $n \in \N$, using $R_n$, $S_n$, $R_n^*$, and $S_n^*$, let us define the correspondence $\mathfrak{R}_n \subset \mathcal{UCP}_{\K}(E_n) \times \mathcal{UCP}_{\K}(E)$ by
	\[\mathfrak R_n = \left\{ ( \phi_n, R_n^* (\phi_n) ) : \phi_n \in \mathcal{UCP}_{\K}(E_n) \right\} \cap \left\{ ( S_n^*(\phi), \phi ) : \phi \in \mathcal{UCP}_{\K}(E) \right\}\]
	with distortion
	\[
	\textup{dis}(\mathfrak R_n) = \sup\left \{ \left| d_{E_n}(\phi_n , \phi_n') - d_E (\phi, \phi') \right| : (\phi_n,  \phi) , ( \phi_n',  \phi')  \in \mathfrak R_n   \right\}.
	\]
	By  \eqref{eqn:dist-GH}, we have 
	\begin{align*}
		d_{GH}(\mathcal{UCP}_{\K}(E_n), \mathcal{UCP}_{\K}(E)) \leq \frac{1}{2} \operatorname{dis } \mathfrak{R}_n.
	\end{align*}
	
	It follows from Proposition \ref{prop:estimate} that $\operatorname{dis } \mathfrak{R}_n$ is bounded by $2(c_n + c_n')$ and hence converges to zero as $c_n$ and $c_n'$ converge to zero. Therefore 
	\[\lim_{n\to\infty}\,d_{GH}(\mathcal{UCP}_{\K}(E_n), \mathcal{UCP}_{\K}(E))=0.\]
	This complete the proof.
	\endproof
	
	A different proof for finite-dimensional $\mathcal K$ can be given using \cite[Proposition 5.19]{Leim} of Leimbach. This was communicated to us by him. He showed a stronger result, viz., an estimate on the complete Gromov-Hausdorff distance (\cite[Definition 3.2]{Kerr_JFA}). The idea of his proof stems from Proposition 2.14 of \cite{KK} which uses the notion of bridge introduced by Rieffel \cite{Rie00}.
	\section{Applications} \label{appl}
\subsection{The Fej\'er--Riesz operator systems}
Let us consider the following spectral triple 
\begin{equation*}
	\left(\E:=C^{\infty}(S^1), \H: = L^2(S^1), D: = - i \frac{d}{d\theta} \right).
\end{equation*}

Consider the Fejer-Riesz operator system $C(S^1)_{(n)}$ (cf. Section \ref{sec:intro}). Since this is an operator subsystem of $E:=C(S^1)$ it is natural to consider the following spectral triple: 

\begin{equation*}
	\left(\E_n:=C(S^1)_{(n)} , \H_n: = L^2(S^1), D_n: = - i \frac{d}{d\theta} \right).
\end{equation*}
That these two spectral triples have the BW property is a consequence of results of section 2.

We will seek a $C^1$-approximate complete order isomorphism $(R_n, S_n)$ for the operator system spectral triples $(\E_n, \H_n, D_n)$ and $(\E, \H, D)$, allowing us to apply Theorem \ref{thm:UCP-GH-conv} to conclude Gromov--Hausdorff convergence of the corresponding spaces of ucp maps.

Let
\begin{align*}
	R_n : C(S^1) &\to C(S^1)_{(n)}\\
	f &\mapsto F_n \ast f 
\end{align*}
where we recall that $F_n = \sum_{|k| \leq n-1} (1-|k|/n ) e^{ik\theta}$ is the Fej\'er kernel (cf. \cite[Chapter 2]{Stein_book}). Also, define
\begin{align*}
	S_n : C(S^1)_{(n)} &\to C(S^1)\\
	f &\mapsto f.
\end{align*}

It follows from Lemma 14 and Lemma 15 of \cite{Suij-JGP} that there exist sequences $\{c_n\}$ and $\{c_n'\}$, both converging to $0$, such that
\begin{align*}
	\| S_n \circ R_n (f) - f \| &\leq c_n \| f \|_1, \\
	\| R_n \circ S_n (h) - h \| &\leq c_n' \| h \|_1.
\end{align*}
From the definition and the properties of the kernel $F_n$, it follows that $R_n$ and $S_n$ are both unital, positive, contractive, and Lipschitz contractive. Since  $C(S^1)$ is a commutative $C^*$-algebra, by \cite[Theorems 3.9 and 3.11]{Pau02}, $R_n$ and $S_n$ are ucp maps. Therefore, $(R_n, S_n)$ form a $C^1$-approximate complete order isomorphism for $(\mathcal{E}_n, \mathcal{H}_n, D_n)$ and $(\mathcal{E}, \mathcal{H}, D)$. Hence, we can apply Theorem \ref{thm:UCP-GH-conv} to these operator systems.

\begin{cor}
	$\mathcal{UCP}_{\K}(C(S^1)_{(n)})$ converges to $\mathcal{UCP}_{\K}(C(S^1))$ in the Gromov-Hausdorff distance.
\end{cor}

\subsection{The Toeplitz operator systems}
Consider the spectral triple $(\mathrm{C}^\infty(S^1), \mathrm{L}^2(S^1), -i\frac{d}{d \theta})$ and the truncated spectral triples 
$$(\mathcal T_n = \{ P_n M_f|_{\mathcal H_n} : f \in \mathrm{C}^\infty(S^1)\}, \mathcal H_n, -iP_n \frac{d}{d \theta}|_{\mathcal H_n}).$$
In the above, $\mathcal H_n$ is the $n$-dimensional space spanned by $\{ e_1, e_2, \ldots , e_n\}$ where $e_k$ is the function $z^k$ on the unit circle and $P_n$ is the orthogonal projection from $ \mathrm{L}^2(S^1)$ onto $\mathcal H_n$. Note that if $f$ has the Fourier series expansion 
$$ f(z) = \sum_{k=-\infty}^{\infty} a_k z^k,$$
then the matrix of $P_n M_f|_{\mathcal H_n}$ in the orthonormal basis $\{ e_1, e_2, \ldots , e_n\}$ is the Toeplitz matrix $(( a_{k-j} ))_{k,j=0}^n$.  Define $R_n$ by sending $f$ in $\mathrm{C}^\infty(S^1)$ to $P_n M_f|_{\mathcal H_n}$. The unit circle has a natural action on $\mathcal T_n$, viz., 
$$\alpha_\theta \left(  (( a_{k-j} ))_{k,j=0}^n \right) = \left(  (( e^{i(k-j) \theta} a_{i-j} ))_{i,j=0}^n \right).$$
Let $\psi$ be the unit vector $(e_1 + e_2 + \cdots + e_n)/\sqrt{n}$. Define $S_n : \mathcal T_n \rightarrow \mathrm{C}(S^1)$ by 
$$(S_n(T))(\theta) =  \langle \psi , \alpha_\theta(T) \psi \rangle.$$
See \cite{Suij-JGP} where it was verified that these are $C^1$ order isomorphisms. Hence we have the following.
\begin{cor}
	$\mathcal{UCP}_{\K}(\mathcal T_n)$ converges to $\mathcal{UCP}_{\K}(C(S^1))$ in the Gromov-Hausdorff distance.
\end{cor}

\subsection{The $d$-dimensional torus}
In this section, $\Omega$ stands for the flat torus $S^1 \times S^1 \times \cdots \times S^1$ ($d$-times). We consider the operator system spectral triple:
\begin{align}\label{eq:d-dim_spt}
	\left(\E := \mathrm{C}^\infty(\Omega), \H := \mathrm{L}^2(\spb(\Omega)), D\right).
\end{align}
Here, $\E$ represents the $*$-algebra of smooth functions on the torus $\Omega$. These functions act by multiplication on a Hilbert space, denoted by $\H$, which consists of $\mathrm{L}^2$-sections of the spinor bundle $\spb(\Omega)$. We identify $\spb(\Omega)$ with the trivial bundle
$\Omega \otimes V$, where $V := \mathbb{C}^d$, and we write $\mathcal{H} := \mathrm{L}^2(\Omega) \otimes V \cong \mathrm{L}^2(\spb(\Omega))$ for the Hilbert space. The Dirac operator $D$ given by 
\[D = -i \sum_{\mu = 1}^{ d } \partial_\mu \otimes \gamma^\mu\] 
acts on the dense subspace of smooth sections of the spinor bundle within $\H$. 
\smallskip

\smallskip
\subsubsection{Spherical truncation}
For an integer $N \geq 0$, let $P_N$ be the orthogonal projection to the subspace of $\mathcal{H}$ spanned by the eigenspinors $e_\lambda$  with $|\lambda| \leq N$. More concretely, we have $P_N \mathcal{H} = \mathrm{span} \{e_n \, : \, n \in \mathbb{Z}^d, \|n\| \leq N\} \otimes V$, with $e_n(x) := e^{i n \cdot x}$, for all $x \in \Omega$. 
\smallskip

Let $\cT^d_N$ be the operator system of $d$-dimensional Toeplitz matrices $T = (t_{k-l})_{k, l \in \bar{\mathrm{B}}_N^\mathbb{Z}}$. Here, $\bar{\mathrm{B}}_N^\mathbb{Z}$ is the set of points in the integer lattice $\mathbb{Z}^d$ that lie inside a ball of radius $N$. The term $t_{k-l}$ represents the inner product $\langle e_k, T e_l \rangle$. Specifically, for $d = 1$, the operator system $\cT^1_N$ corresponds to the set of Toeplitz matrices of size $(2 N +1) \times (2 N +1)$, which were studied in detail in \cite{CS20} and \cite{Hek22}.

The spectral projection $P_N$ gives rise to the following \emph{operator system spectral triple}:
\begin{align}\label{eq:sph_tr}
	\left(\E_N:= \cT^d_N, \H_N:=P_N \mathcal{H}, D_N:=P_N D |_{\H_N} \right),
\end{align}
where $\H$ and $ D$ are given by \eqref{eq:d-dim_spt}. {\em We call the above truncation \eqref{eq:sph_tr} the spherical truncation of the \eqref{eq:d-dim_spt}.}  Note that the truncations for the torus are different from the truncation for the unit circle discussed above. We recall the following from \cite{LeSu24}.

Define 
\begin{align*}
	R_N:&~\E\to\E_N \\
	&f\mapsto P_NfP_N,\\[5pt]
	S_N:&~\E_N\to \E \\
	&T\mapsto \frac{1}{\mathcal{N}_\mathrm{B}(N)} \mathrm{Tr} \left( \ket{\psi} \bra{\psi} \alpha(T) \right).
\end{align*}
Here, $\alpha$ is the $\Omega$-action on $\E_N$ given by 
\begin{align*}
	\alpha_\theta(T) = \left( t_{k-l} e^{i (k-l) \cdot \theta}\right)_{k,l \in \bar{\mathrm{B}}_N^\mathbb{Z}},
\end{align*}
where $\theta = (\theta_1, \theta_2, \cdots, \theta_d)$ such that $(e^{i\theta_1}, e^{i\theta_2}, \cdots, e^{i\theta_d}) \in \Omega$, the vector $\ket{\psi}$ is given by $\ket{\psi} = \sum_{n \in \bar{\mathrm{B}}_N^\mathbb{Z}} e_n$, and $\mathcal{N}_\mathrm{B}(N) := | \bar{\mathrm{B}}_N \cap \mathbb{Z}^d | $ is the number of $\mathbb{Z}^d$-lattice points in the closed ball of radius $N$. By Lemma \ref{lem:BW=d} and Lemma \ref{lem:BW=d:1}, $(\E_N, \H_N, D_N)$ and $(\E, \H, D)$ have the {\it BW}-property.

It is shown in the proof of \cite[Theorem 3.10]{LeSu24} that the pair $(R_N, S_N)$ form a $C^1$-approximate order isomorphism for the spectral triples $(\E_N, \H_N, D_N)$ and $(\E, \H, D)$. By \cite[Theorems 3.9 and 3.11]{Pau02}, $R_N$ and $S_N$ are ucp maps, and hence the pair $(R_N, S_N)$ form a $C^1$-approximate complete order isomorphism. Therefore, applying Theorem \ref{thm:UCP-GH-conv}, we have the following convergnce result.

\begin{cor}
	Let $d\in \N$. Then $\mathcal{UCP}_{\K}(\cT_N^d)$ converges to $\mathcal{UCP}_{\K}(C(\Omega))$ in the Gromov-Hausdorff distance.
\end{cor}

\subsubsection{Polyhedral truncation}
This section is inspired by \cite{Trav}.

\begin{defn}
	Let $r,d \in \mathbb{N}$. A \textit{polyhedron} $\tilde{\Delta}(r)$ in $\mathbb{R}^d$ is the convex hull of a finite set $W=\{ \bm a_1, \ldots, \bm a_r \} \subset \mathbb{Z}^d$, where $W$ is minimal, meaning that it coincides with the set of vertices of $\tilde{\Delta}(r)$. The restriction $\Delta(r) = \tilde{\Delta}(r) \cap \mathbb{Z}^d$ is called a \textit{summability polyhedron} if $\tilde{\Delta}(r)$ contains the origin strictly in its interior. For any $N \in \mathbb{N}$, we denote by $\tilde{\Delta}_N(r)$ the convex hull of the set $\{ N \bm a_1, \ldots, N \bm a_r \}$ and define $\Delta_N(r) = \tilde{\Delta}_N(r) \cap \mathbb{Z}^d$. When the polyhedron is fixed, that is, $r$ and $W$ are fixed, we simply write $\Delta$ and $\Delta_N$ instead of $\Delta(r)$ and $\Delta_N(r)$, respectively. We call $\Delta_N$  the $N$-dilation of $\Delta$.
	
\end{defn}

\smallskip
Let $\Delta$ be a summability polyhedron and let  $\Delta_N$ be the $N$-dilation of $\Delta$ for an integer $N \ge 1$. Let $Q_N$ be the orthogonal projection to the subspace of $\mathcal{H}$ spanned by the eigenspinors $e_\lambda$ with $\lambda\in\Delta_N$, i.e.,  $Q_N \mathcal{H} = \mathrm{span} \{e_\lambda \, : \lambda\in\Delta_N\} \otimes V$, with $e_n(x) := e^{i n \cdot x}$, for all $x \in \Omega$. 
\smallskip

Let $\cT^d_N(PH)$ (PH stands for polyhedron) be the operator system of $d$-dimensional polyhedral Toeplitz matrices $T = (t_{k-l})_{k, l \in \Delta_N}$. This is a $| \Delta_N | \times | \Delta_N |$ matrix with entries indexed by the elements of $\Delta_N$. The term $t_{k-l}$ represents the inner product $\langle e_k, T e_l \rangle$.

The spectral projection $Q_N$ defines the \emph{operator system spectral triple}:
\begin{align}\label{eq:phdrl_tr}
	\left(\E_N:= \cT^d_N(PH), \H_N:=Q_N \mathcal{H}, D_N:=Q_N D |_{\H_N} \right),
\end{align}
where $\H$ and $ D$ are given by \eqref{eq:d-dim_spt}. {\em We call the above truncation \eqref{eq:phdrl_tr} a polyhedral truncation of \eqref{eq:d-dim_spt}.}  Define 
\begin{align*}
	\rho_N:&~\E\to\E_N \\
	&f\mapsto Q_NfQ_N,\\[5pt]
	\sigma_N:&~\E_N\to \E \\
	&T\mapsto \frac{1}{| \Delta_N |} \mathrm{Tr} \left( \ket{\psi} \bra{\psi} \alpha(T) \right).
\end{align*}
Here, $\alpha$ is the $\Omega$-action on $\E_N$ given by 
\begin{align*}
	\alpha_\theta(T) = \left( t_{k-l} e^{i (k-l) \cdot \theta}\right)_{k,l \in \Delta_N},
\end{align*}
where $\theta = (\theta_1, \theta_2, \cdots, \theta_d)$ is such that $(e^{i\theta_1}, e^{i\theta_2}, \cdots, e^{i\theta_d}) \in \Omega$ and the vector $\ket{\psi}$ is given by $\ket{\psi} = \sum_{n \in \Delta_N} e_n$. By Lemma \ref{lem:BW=d} and Lemma \ref{lem:BW=d:1}, $(\E_N, \H_N, D_N)$ and $(\E, \H, D)$ have the {\it BW}-property.

Using the approach of \cite{LeSu24}, we show below that the pair $(\rho_N, \sigma_N)$ forms a $C^1$-approximate order isomorphism for the spectral triples $(\E_N, \H_N, D_N)$ and $(\E, \H, D)$.

\begin{lem}
	The maps $\rho_N$ and $\sigma_N$ are unital, completely positive and $\mathrm{C}^1$-contractive.
\end{lem}
\begin{proof}
	The proof is similar to \cite[Lemma 3.1 and Lemma 3.2]{LeSu24}, along with an application of \cite[Theorems 3.9 and 3.11]{Pau02} as $C^\infty(\Omega)$ is commutative.
\end{proof}
\begin{prop}
	The pair $(\rho_N, \sigma_N)$ forms a $C^1$-approximate order isomorphism for the spectral triples $(\mathcal{E}_N, \mathcal{H}_N, D_N)$ and $(\mathcal{E}, \mathcal{H}, D)$ given by \eqref{eq:phdrl_tr} and \eqref{eq:d-dim_spt}, respectively.
\end{prop}
\begin{proof}
	By a computation similar to \cite[Lemma 3.3]{LeSu24}, we have 
	\begin{align*}
		&\sigma_N \circ \rho_N (f) (x)=\sum_{n \in \Delta_N-\Delta_N} \frac{| \Delta_N\cap(\Delta_N+n)|}{| \Delta_N |} \widehat{f}(n) e_n(x)\\
		&\rho_{N} \circ \sigma_N (T)= \left( \frac{ | \Delta_N\cap(\Delta_N+m-n)|}{| \Delta_N |} t_{m-n} \right)_{m,n \in\Delta_N}
	\end{align*}
For each $N\in\N$, define $K_{2N}$ on $\Omega$ by
\begin{align*}
	K_{2N}(x)=&\sum_{n,m \in \Delta_N} \frac{1}{| \Delta_N |}  e_{n-m}(x)\\
	=&\sum_{n \in \Delta_N-\Delta_N} \frac{| \Delta_N\cap(\Delta_N+n)|}{| \Delta_N |} e_n(x).
\end{align*}
Then it is easy to check that 
\[\sigma_N \circ \rho_N (f) (x)=(K_{2N}*f)(x).\]
A careful examination of the proof of \cite[Theorem 3.10]{LeSu24} suggests that, to demonstrate that $(\rho_N, \sigma_N)$ form a $C^1$-approximate order isomorphism, it suffices to show that the family $\{K_{2N}\}_{N \in \mathbb{N}}$ serves as a good kernel. In \cite{Trav}, $K_{2N}$ is called the polyhedron Fejér kernel. From \cite[Theorems 2, 3]{Trav}, it follows that:
\begin{enumerate}
	\item $K_{2N} \geq 0$, 
	\vspace{.05in} 
	\item $ \int_{\Omega} K_{2N}(x) \, dx = 1$. 
	\vspace{.05in}
	\item $\hat{K}_{2N}\to 1$ pointwise on $\Z^d$.
\end{enumerate}  
Next, we need that the mass of $K_{2N}$ is concentrated near the origin. But this can be proved by mimicking the proof of \cite[Lemma 3.8]{LeSu24} and we omit the details. Hence $\{K_{2N}\}$ is good kernel. This concludes the proof.
\end{proof}

\begin{cor}
	Let $d\in \N$. Then $\mathcal{UCP}_{\K}(\cT_N^d(PH))$ converges to $\mathcal{UCP}_{\K}(C(\Omega))$ in the Gromov-Hausdorff distance.
\end{cor}

\begin{rem}
  As one more example of Gromov-Hausdorff convergence of CP maps, consider $V_n$, the $n$-dimensional irreducible representation of $SU(2)$. Let $S^2$ be the $2$-dimensional unit sphere. Then $\mathcal{UCP}_\K (L(V_n))$ converges to $\mathcal{UCP}_\K (C(S^2))$. This is a consequence of non-trivial results of van Suijlekom in section 3.3 of \cite{Suij-JGP} and our Theorem \ref{thm:UCP-GH-conv}. van Suijlekom produces the $C^1$-approximate order isomorphism and hence Theorem \ref{thm:UCP-GH-conv} can be applied.
\end{rem}

 \textsf{Acknowledgement: The authors are thankful to Walter van Suijlekom, whose India visit under the Jubilee Chair Professorship of the Indian of Academy of Sciences was pivotal for this paper, as well as M. Leimbach for his comments. The first author is supported by a J C Bose Fellowship JCB/2021/000041 of SERB. The second author is supported by the Prime Minister’s Research Fellowship. The third author is supported by NBHM post-doctoral fellowship.}}
 
 \section*{Declarations}
 
 On behalf of all authors, the corresponding author states that there is no conflict of interest.
 
 This manuscript has no associated data.


\begin{thebibliography}{99}
	
	\bibitem{BBI01}
	D.~Burago, Y.~Burago, and S.~Ivanov.
	\newblock {\em A course in metric geometry}, volume~33 of {\em Graduate Studies
		in Mathematics}.
	\newblock American Mathematical Society, Providence, RI, 2001.
	
	
	\bibitem{CS20}
	A.~Connes and W.~van Suijlekom.
	\newblock Spectral truncations in noncommutative geometry and operator systems.
	\newblock {\em Commun. Math. Phys.} 383 no. 3 (2021), 2021--2067.
	
	
	\bibitem{Far}
	D. ~Farenick.
	\newblock The operator system of Toeplitz matrices.
	\newblock {\em Trans. Amer. Math. Soc. Ser. B} 8 (2021), 999–1023.
	
	\bibitem{FLP}
	C. ~Farsi, F. ~Latrémolière, and J. ~Packer.
	\newblock Convergence of inductive sequences of spectral triples for the spectral propinquity.
	\newblock {\em Adv. Math.} 437 (2024), Paper No. 109442, 59 pp.
	
	\bibitem{HHIK}
	Y. ~Hashimoto, A. ~Hafid, M. ~Ikeda, H. ~Kadri.
	\newblock Spectral Truncation Kernels: Noncommutativity in $C^*$-algebraic Kernel Machines.
	\newblock Preprint (2024), \url{https://doi.org/10.48550/arXiv.2405.17823}.
	
	
	\bibitem{Hek22}
	E. Hekkelman.
	\newblock Truncated geometry on the circle.
	\newblock {\em Lett. Math. Phys.} 112 no. 2 (2022) Paper No. 20, 19pp.
	
	\bibitem{KK}
	J. Kaad, D. Kyed.
	\newblock The quantum metric structure of quantum $SU(2)$.
	\newblock  Memoirs of the EMS, 2025. 
	
	\bibitem{Kerr_JFA}
	D. Kerr.
	\newblock Matricial quantum Gromov-Hausdorff distance.
	\newblock {\em J. Funct. Anal.} 205 (2003), 132--167.

\bibitem{KL}
D. ~Kerr, H. Li.
\newblock On gromov-Haudorff convergence for operator metric spaces.
\newblock {J. Operator Thery} 62 (2009), 83--109. 

\bibitem{Leim}
M. Leimbach.
\newblock Convergence of Peter--Weyl Truncations of Compact Quantum Groups.
\newblock Preprint (2024), \url{https://doi.org/10.48550/arXiv.2409.16698}.
	
	
	\bibitem{LeSu24}
	M.~Leimbach, W.~van Suijlekom.
	\newblock Gromov-{H}ausdorff convergence of spectral truncations for tori.
	\newblock {\em Adv. Math.} 439 (2024) Paper No. 109496, 26pp.
	
	\bibitem{Pau02}
	V.~Paulsen.
	\newblock {\em Completely bounded maps and operator algebras}, volume~78 of
	{\em Cambridge Studies in Advanced Mathematics}.
	\newblock Cambridge University Press, Cambridge, 2002.
	
	
	\bibitem{Rie00}
	M.~A. Rieffel.
	\newblock Gromov-{H}ausdorff distance for quantum metric spaces.
	\newblock {\em Mem. Amer. Math. Soc.} 168 (2004)  1--65.
	
	
	\bibitem{Rie2023}
	M.~A. Rieffel.
	\newblock  Convergence of Fourier truncations for compact quantum groups and finitely generated groups.
	\newblock {\em J. Geom. Phys.} 192 (2023), Paper No. 104921, 13 pp.
	
	\bibitem{Stein_book}
	E. Stein, R. Shakarchi. 
	\newblock {\em Fourier Analysis: An Introduction.} 
	\newblock Princeton University
	Press, 2003.
	
	
	\bibitem{Suij-JGP}
	W.~van Suijlekom.
	\newblock Gromov-{H}ausdorff convergence of state spaces for spectral truncations.
	\newblock {\em J. Geom. Phys.} 162 no. 2 (2021)  Paper No. 104075, 11pp.
	
	\bibitem{Trav}
	G. Travaglini.
	\newblock Fej\'er kernels for {F}ourier series on {${\bf T}^n$} and on compact {L}ie groups.
	\newblock {\em Math. Z.} 216 (1994) 265--281
\end{thebibliography}
\end{document}